\documentclass[a4paper,12pt,reqno]{amsart}
\usepackage[top=29truemm,bottom=24truemm,left=23truemm,right=22truemm]{geometry}
\usepackage{amsfonts,amssymb,amsmath,amsthm,graphicx,color,bm,ulem,booktabs}
\usepackage{latexsym,subfigure,ascmac,lscape,cite}
\definecolor{ultramarine}{rgb}{0.07, 0.04, 0.56}
\definecolor{cadmiumgreen}{rgb}{0.0, 0.42, 0.24}
\definecolor{indigo(dye)}{rgb}{0.0, 0.25, 0.42}
\usepackage[linktocpage=true]{hyperref}
\hypersetup{
colorlinks=true,
citecolor=ultramarine,
linkcolor=cadmiumgreen,
urlcolor=indigo(dye),
}

\theoremstyle{plain}
\newtheorem{theorem}{Theorem} 
\newtheorem{corollary}{Corollary}[theorem]
\newcommand{\f}[2]{\dfrac{#1}{#2}}  
\newcommand{\mk}[1]{\left( #1 \right)}  
\newcommand{\kk}[1]{\left[ #1 \right]}  
  
\newcommand{\be}{\begin{equation}}  
\newcommand{\ee}{\end{equation}}
\newcommand{\bem}{\begin{matrix}}  
\newcommand{\eem}{\end{matrix}}
\newcommand{\Func}[5]{{}_{#1}F_{#2}\mk{ \bem #3 \\ #4 \eem ;#5}}
\newcommand{\ba}{{\bf a}}
\newcommand{\bb}{{\bf b}}
\newcommand{\sZ}{\mathbb{Z}}
\newcommand{\sC}{\mathbb{C}}
\newcommand{\sN}{\mathbb{N}}

\begin{document}

\title{Differentiation identities for hypergeometric functions}

\author{Hayato Motohashi}
\address{Division of Liberal Arts, Kogakuin University, 2665-1 Nakano-machi, Hachioji, Tokyo, 192-0015, Japan}
\email{motohashi@cc.kogakuin.ac.jp}


\begin{abstract}
It is well-known that differentiation of hypergeometric function multiplied by a certain power function yields another hypergeometric function with a different set of parameters.
Such differentiation identities for hypergeometric functions have been used widely in various fields of applied mathematics and natural sciences.
In this expository note, we provide a simple proof of the differentiation identities, which is based only on the definition of the coefficients for the power series expansion of the hypergeometric functions.
\end{abstract}

\maketitle

\section{Introduction}
\label{sec:intro}

The theory of hypergeometric functions is fundamental in 
mathematical physics
since they contain as special cases almost all the special functions commonly used for centuries.
Among the family of hypergeometric functions, the most well-known ones are Gauss hypergeometric function and Kummer confluent hypergeometric function, 
and their natural extension is known as 
generalized hypergeometric function~\cite{Bailey,Slater}.
A lot of efforts have been made to deepen understanding of the hypergeometric functions, which have accumulated hundreds of useful formulae~\cite{Bateman:100233,mabramowitz64:handbook,Luke,NIST:DLMF}, 
such as relation to other special functions, series expansions, differential properties, contiguous relations, variable transformations, integral representations etc.
These formulae clarify the nature of the hypergeometric functions from various aspects, and have been widely used as powerful tools in diverse research fields~\cite{Seaborn}.

Among such formulae, let us focus on the differentiation identity for the hypergeometric functions.
Consider the Gauss hypergeometric function as an example, which is defined by an infinite series
\be \label{2F1def} \Func{2}{1}{a,b}{c}{z} = \sum_{k=0}^\infty \f{1}{k!}\f{(a)_k(b)_k}{(c)_k} z^k , \ee
where $(a)_k$ is the Pochhammer symbol defined by 
\be \label{Pdef} (a)_k=\f{\Gamma(a+k)}{\Gamma(a)}=
\begin{cases}
      1, & (k=0), \\
      a(a+1)\cdots (a+k-1), & (k>0).
\end{cases}
\ee
Here, it is assumed that $c$ is not a nonpositive integer to avoid the singular behavior of the coefficients in \eqref{2F1def}.
The series converges on $|z|<1$.
Several differentiation identities are well-known in the literature.
One of the famous identities is as follows (see, e.g., \cite[Eq.~(15.2.4)]{mabramowitz64:handbook}):
\be \label{Th1-4-old} \f{d^n}{dz^n}\kk{z^{c-1}\Func{2}{1}{a,b}{c}{z}} = (c-n)_n z^{c-n-1} \Func{2}{1}{a,b}{c-n}{z}. \ee
While it is typically not clearly mentioned in the literature, this identity does not hold if $c\in \sZ_{< n+1}$.
For $c\in \sZ_{< n+1}$, the right-hand side is singular.
For $c\in \sZ_{\leq 0}$, the left-hand side is also singular.
However, for $c=1,2,\cdots, n$, the left-hand side is not singular, despite that the right-hand side is singular.
It is thus natural to ask whether there exists a different kind of identity for the exceptional case $c=1,2,\cdots, n$ to rewrite the left-hand side in terms of a hypergeometric function with another set of parameters.
While most books do not make any mention of such exceptional cases, there indeed exists different forms of the differentiation identities for these cases, some of which can be found in \cite[page 44]{Luke}.

\begin{table}[t]
 \begin{center}
   \caption{Numerical evaluation of the left-hand side ($f_{L}$) and the first ($f_{R1}$) and second ($f_{R2}$) lines of the right-hand side of \eqref{Th1-4-sp} for several values of $c$ with a specific parameter set $n=4$, $a=1/2$, $b=2/3$, and $z=1/3$.}
   \begin{tabular}{cccc}
      \toprule
      $c$ & $f_L$ & $f_{R1}$ & $f_{R2}$ \\ \midrule
      $1$ & $16.2802578209098$ & - & $16.2802578209098$ \\ 
      $2$ & $3.39340187542396$ & - & $3.39340187542396$ \\
      $3$ & $2.04681438609744$ & - & $2.04681438609744$ \\
      $4$ & $3.31081155003091$ & - & $3.31081155003091$ \\
      $5$ & $27.4105535888826$ & $27.4105535888826$ & $27.4105535888826$ \\
      $6$ & $42.6040520193532$ & $42.6040520193532$ & - \\
      $7$ & $41.6637846070299$ & $41.6637846070299$ & - \\
      \bottomrule
   \end{tabular}
   \label{table:Th1-4}
 \end{center}
\end{table}

\begin{figure}[t]
 \begin{center}
   \includegraphics[clip,width=0.76\columnwidth]{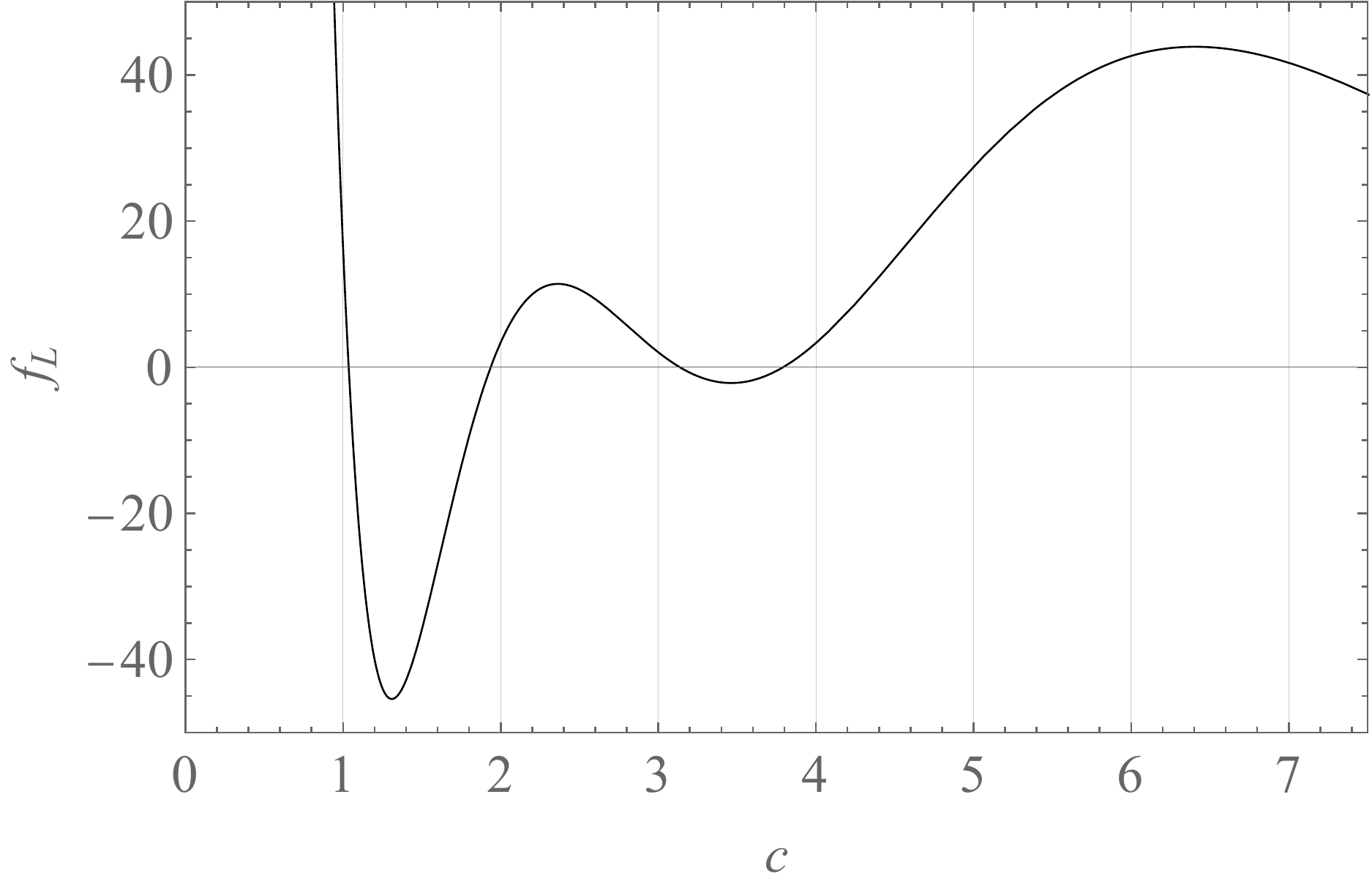}
   \includegraphics[clip,width=0.76\columnwidth]{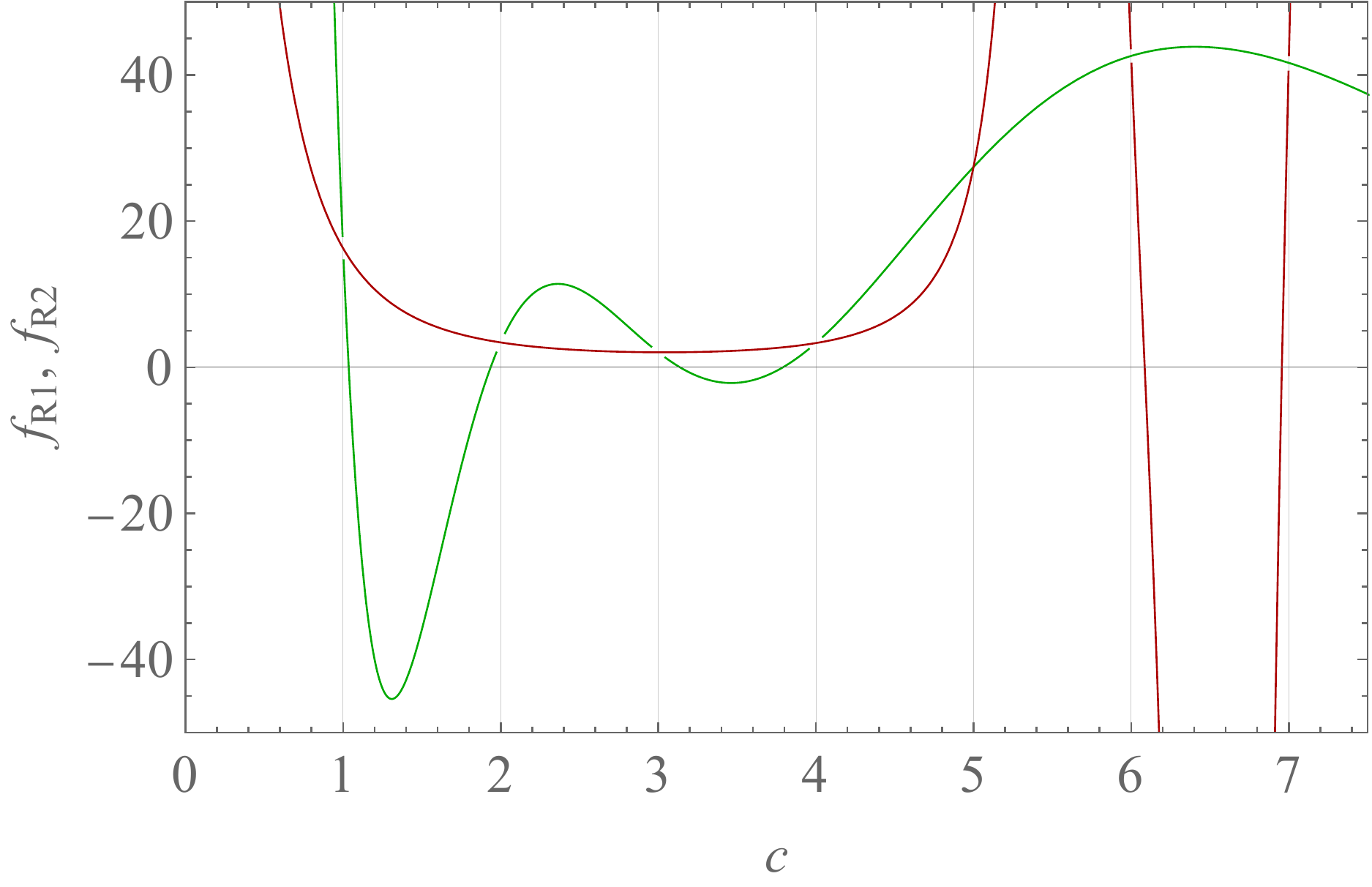}
   \caption{The left-hand side ($f_{L}$, black) and the first ($f_{R1}$, green) and second ($f_{R2}$, red) lines of the right-hand side of \eqref{Th1-4-sp} for a specific parameter set $n=4$, $a=1/2$, $b=2/3$, and $z=1/3$.}
   \label{fig:Th1-4}
 \end{center}
\end{figure}

By supplementing the exceptional case, the identity~\eqref{Th1-4-old} is improved as 
(see \eqref{Th1-4} below)
{\fontsize{10pt}{10pt}\selectfont 
\begin{align}
&\f{d^n}{dz^n}\kk{z^{c-1} \Func{2}{1}{a,b}{c}{z} } \notag\\
&= 
\begin{cases}
      (c-n)_n z^{c - n - 1} \Func{2}{1}{a,b}{c-n}{z}, & (c\not\in \sZ_{< n+1}) , \\
      \f{n!}{(n-c+1)!}\f{(a)_{n-c+1}(b)_{n-c+1}}{(c)_{n-c+1}} \Func{2}{1}{a+n-c+1,b+n-c+1}{n-c+2}{z} , & (c=1,2,\cdots,n+1) .
\end{cases} \label{Th1-4-sp}
\end{align}}
\noindent
Note that the first and second lines of the right-hand side are identical for $c=n+1$.
In Table~\ref{table:Th1-4}, we present example numerical values of the left-hand side ($f_{L}$) and the first ($f_{R1}$) and second ($f_{R2}$) lines of the right-hand side of \eqref{Th1-4-sp} for several values of $c$, while we fix the other parameters to a specific parameter set $n=4$, $a=1/2$, $b=2/3$, and $z=1/3$.
They match each other for the range of the validity of the identity~\eqref{Th1-4-sp}.
As explained above, $f_{L}$ is regular but $f_{R1}$ is singular for $c=1,2,3,4$, and $f_{R2}$ precisely complements the identity for the exceptional case.
On the other hand, for $c=6,7,\cdots$, $f_{R1}$ is regular but $f_{R2}$ is singular since the lower parameter $n-c+2$ in $f_{R2}$ is a nonpositive integer. 
In Fig.~\ref{fig:Th1-4}, we depict the left-hand side ($f_{L}$, black) and the first ($f_{R1}$, green) and second ($f_{R2}$, red) lines of the right-hand side of \eqref{Th1-4-sp} as a function of the parameter $c$, with the same parameter set.
The $f_{R1}$ and $f_{R2}$ intersect each other at $c=5$.
For $c=1,2,3,4$ and $c=6,7,\cdots$, one of $f_{R1}$ and $f_{R2}$ is singular, but we see an interesting behavior where the regular one complements the singular one.
Therefore, the difference $f_{R1}-f_{R2}$ has a peculiar periodic behavior as a function of $c$. 

Of course, one can derive the formula for the case $c=1,2,\cdots, n$ by applying the analytic continuation to the right-hand side of \eqref{Th1-4-old}.
It is then possible to enlarge the validity of the formula and to obtain interesting pairs of functions demonstrated above.
However, such a derivation is based on the assumption that one already knows the original formula~\eqref{Th1-4-old}.
Thus, this strategy may not work for more general functions, if no differentiation identities are known for those functions. 

The aim of this expository note is to provide a simple derivation of the differentiation identities for the generalized hypergeometric function, treating the well-known cases such as \eqref{Th1-4-old} and lesser-known exceptional cases on an equal footing.
The advantage of this derivation is that, 
unlike the approach relying on the analytic continuation mentioned above, it is based only on the definition of the coefficients for the power series expansion of the hypergeometric functions and does not require a priori knowledge of the differentiation identities.
The strategy presented here may work for more general functions~\cite{Motohashi:prep}, while the coefficients for the power series expansion are in general not so simple.

In the rest of the paper, we derive differentiation identities for the generalized hypergeometric function, including the hypergeometric function and confluent hypergeometric function as special cases, in a concise and systematic way to exhaust all the possible identities, some of which are omitted in the literature.

\section{Prelude}
\label{sec:pre}

In this section we summarize preliminary definitions and notations.
Generalized hypergeometric function is defined formally by an infinite series 
\be \label{pFqdef} \Func{p}{q}{a_1,\cdots,a_p}{b_1,\cdots,b_q}{z} = \sum_{k=0}^\infty c_k z^k , \ee
with coefficients
\be \label{ckdef} c_k= \f{1}{k!}\f{(a_1)_k\cdots (a_p)_k}{(b_1)_k\cdots (b_q)_k}, \ee
with the Pochhammer symbol $(a)_k$ defined in \eqref{Pdef}.
Here, $p,q \in \sN$, and other parameters and the variable $z$ are complex numbers in general.

The case with $(p,q)=(2,1)$ yields the hypergeometric function, whereas the case with $(p,q)=(1,1)$ yields the confluent hypergeometric function.
For these cases, the notations $F(a,b;c;z)=\Func{2}{1}{a,b}{c}{z}$ and $M(a;c;z)=\Func{1}{1}{a}{c}{z}$ are also commonly used in the literature.
These two hypergeometric functions are especially important and have been extensively studied since they are exact solutions for second-order differential equations known as (confluent) hypergeometric equations, which encompass a broad class of differential equations appearing in mathematical physics.

For simplicity, we use the notation $\ba=(a_1,\cdots,a_p)$ and $\bb=(b_1,\cdots,b_q)$ and express the generalized hypergeometric function~\eqref{pFqdef} as $\Func{p}{q}{\ba}{\bb}{z}$.
We also adopt the notation $\ba+c=(a_1+c,\cdots,a_p+c)$ and $\bb+c=(b_1+c,\cdots,b_q+c)$.
Further, we write the product of the Pochhammer symbols as $(\ba)_k=(a_1)_k\cdots (a_p)_k$ and $(\bb)_k=(b_1)_k\cdots (b_q)_k$.
The coefficient~\eqref{ckdef} is then written as $c_k=\f{1}{k!}\f{(\ba)_k}{(\bb)_k}$.
By definition, the value of $\Func{p}{q}{\ba}{\bb}{z}$ remains the same under any permutation of the arguments of $\ba$, and any permutation of the arguments of $\bb$.
Therefore, for the following, when we focus on an argument of $\ba$ or $\bb$, we always take the first argument $a_1$ or $b_1$ without loss of generality, and write the rest arguments as $\hat \ba=(a_2,a_3,\cdots,a_p)$ or $\hat \bb=(b_2,b_3,\cdots,b_q)$.
Note that $\hat \ba$ is empty for $p=1$, and $\hat \bb$ is empty for $q=1$.

The radius of convergence of the generalized hypergeometric series~\eqref{pFqdef} is as follows~\cite{Bailey,Slater}:
If $p<q+1$, the series converges for all values of $z$ on the complex plane, provided that the coefficients~\eqref{ckdef} are not singular.
If $p>q+1$, the series does not converge except for $z=0$, and hence the function is defined only when the series terminates.
If $p=q+1$, the series converges when $|z|<1$.
Whether it converges for the case $p=q+1$ when $|z|=1$ is more subtle.
It converges when $z=1$ if 
\be \Re\mk{ \sum_{j=1}^q b_j -\sum_{j=1}^p a_j } >0, \ee 
and when $z=-1$ if 
\be \Re\mk{ \sum_{j=1}^q b_j -\sum_{j=1}^p a_j } + 1 >0. \ee

If none of the upper parameters $\ba$ is a nonpositive integer, we 
assume that none of the lower parameters $\bb$ is a nonpositive integer to avoid the possibility of singular coefficients.
On the other hand, if one or more of the upper parameters $\ba$ are nonpositive integer, the series~\eqref{pFqdef} terminates and $\Func{p}{q}{\ba}{\bb}{z}$ becomes a polynomial.
For instance, if $a_1=-m$ such that $m \in \sZ_{\geq 0}$, the coefficient $c_k$ vanishes for $k\geq m+1$, and hence $\Func{p}{q}{\ba}{\bb}{z}$ is at most $m$-th order polynomial.
Namely, we can write down the generalized hypergeometric function as a finite series as
\be \label{negative} \Func{p}{q}{-m,\hat \ba}{\bb}{z} = \sum_{k=0}^m c_k z^k. \ee
This expression is also valid even if $b_1=-m-\ell$ such that $\ell \in \sZ_{\geq 0}$ provided that we define 
\be \label{inter} \Func{p}{q}{-m,\hat \ba}{-m-\ell,\hat \bb}{z} = \lim_{b_1\to -m-\ell}\mk{ \lim_{a_1\to-m} \Func{p}{q}{a_1,\hat \ba}{b_1,\hat \bb}{z} } . \ee
This definition is a natural generalization of the one for the Gauss hypergeometric function, which is commonly used in the literature, see, e.g., 
\cite[Eq.~(15.2.5)]{NIST:DLMF} or \cite[page 39]{Magnus}.
Throughout the present paper, we adopt the definition~\eqref{inter}.
In other words, for the case of \eqref{negative}, we relax the assumption on the lower parameters $\bb$ as follows: None of the lower parameters $\bb$ is a negative integer $-k$ such that $k \in \sZ_{<m}$.

By definition, it holds that
\be \label{cancel} \Func{p+1}{q+1}{c,\ba}{d,\bb}{z} = \Func{p}{q}{\ba}{\bb}{z}, \ee
if $c=d$.
For brevity, for the following we call this reduction procedure ``cancellation'' between $c$ and $d$.

\section{Main Theorem}

\begin{landscape}
\begin{theorem}\label{thm1}
For $n\in \sN$ and $r\in \sC$, the following differentiation identities for the generalized hypergeometric function ${}_{p}F_{q}$ hold:
\begin{align}
&\f{d^n}{dz^n}\kk{z^r \Func{p}{q}{\ba}{\bb}{z} } \notag\\
&= 
\begin{cases}
      (r-n+1)_nz^{r-n} \Func{p+1}{q+1}{r+1,\ba}{r-n+1,\bb}{z}, & (r\not\in \sZ_{< n}) , \\
      \f{n!}{(n-r)!}\f{(\ba)_{n-r}}{(\bb)_{n-r}} \Func{p+1}{q+1}{n+1,\ba+n-r}{n-r+1,\bb+n-r}{z} , & (r=0,1,\cdots,n) , \\
      (r-n+1)_nz^{r-n} \Func{p+1}{q+1}{r+1,\ba}{r-n+1,\bb}{z} + \f{n!}{(n-r)!}\f{(\ba)_{n-r}}{(\bb)_{n-r}} \Func{p+1}{q+1}{n+1,\ba+n-r}{n-r+1,\bb+n-r}{z} , & (r\in \sZ_{<0}),
\end{cases} \label{Th1-1} \\
&\f{d^n}{dz^n} \Func{p}{q}{\ba}{\bb}{z} = \f{(\ba)_{n}}{(\bb)_{n}} \Func{p}{q}{\ba+n}{\bb+n}{z} , \label{Th1-2} \\
&\f{d^n}{dz^n}\kk{z^{a_1+n-1} \Func{p}{q}{\ba}{\bb}{z} } = (a_1)_n z^{a_1-1} \Func{p}{q}{a_1+n,\hat\ba}{\bb}{z} , \label{Th1-3} \\
&\f{d^n}{dz^n}\kk{z^{b_1-1} \Func{p}{q}{\ba}{\bb}{z} } \notag\\
&= 
\begin{cases}
      (b_1-n)_n z^{b_1 - n - 1} \Func{p}{q}{\ba}{b_1-n,\hat\bb}{z}, & (b_1\not\in \sZ_{< n+1}) , \\
      \f{n!}{(n-b_1+1)!}\f{(\ba)_{n-b_1+1}}{(\bb)_{n-b_1+1}} \Func{p}{q}{\ba+n-b_1+1}{n-b_1+2,\hat\bb+n-b_1+1}{z} , & (b_1=1,2,\cdots,n+1) ,
\end{cases} \label{Th1-4} \\
&\f{d^n}{dz^n}\kk{z^r \Func{p}{q}{1,\hat\ba}{\bb}{z} } \notag\\
&= 
\begin{cases}
\f{n!(\hat \ba)_{n-r}}{(\bb)_{n-r}} \Func{p}{q}{n+1,\hat\ba+n-r}{\bb+n-r}{z}, & (r=0,1,\cdots,n) , \\
\f{n!(\hat \ba)_{n-r}}{(\bb)_{n-r}} \Func{p}{q}{n+1,\hat\ba+n-r}{\bb+n-r}{z}
+ (r-n+1)_n z^{r-n} \Func{p+1}{q+1}{r+1,1,\hat\ba}{r-n+1, \bb}{z} , & (r\in \sZ_{<0}) . 
\end{cases} \label{Th1-5}
\end{align}
\end{theorem}
\end{landscape}

We can easily check that the identities in Theorem~\ref{thm1} (and Corollaries~\ref{coy1} and \ref{coy2} below) hold numerically along the same line of Table~\ref{table:Th1-4} and Fig.~\ref{fig:Th1-4}, while we do not repeat them here.

Note that, similarly to the case $c=n+1$ for \eqref{Th1-4-sp}, there are some cases where the formulae are identical; 
The first and second lines of the right-hand side of \eqref{Th1-1} with $r=n$, 
\eqref{Th1-4} with $b_1=n+1$, and 
\eqref{Th1-5} with $r=0$ are identical, respectively.

Note also that, while different kinds of the generalized hypergeometric functions, ${}_pF_q$ and ${}_{p+1}F_{q+1}$, appear on both sides of the above identities, as expected, the term-by-term differentiation of a convergent power series has the same radius of convergence as the original power series.  Indeed, both radii of convergence of ${}_pF_q$ and ${}_{p+1}F_{q+1}$ are determined by the value of $p-q=(p+1)-(q+1)$, as explained in \S\ref{sec:pre}.

While this type of differentiation identities can be found in the literature, some of them are not written down explicitly.
Specifically, the first line of \eqref{Th1-1}, \eqref{Th1-2}, \eqref{Th1-3}, and the first line of \eqref{Th1-4} correspond to, e.g., Eqs.~(16.3.1)--(16.3.4) in \cite{NIST:DLMF}.
The second line of \eqref{Th1-1} corresponds to \cite[page 44, Eq.~(3)]{Luke}.
However, to the best of my knowledge, it seems that there is no reference that explicitly mentions the third line of \eqref{Th1-1}, the second line of \eqref{Th1-4}, and \eqref{Th1-5}.
The situation is the same for the confluent hypergeometric function with $(p,q)=(1,1)$ and the hypergeometric function with $(p,q)=(2,1)$.
For these cases, \eqref{Th1-2}, \eqref{Th1-3}, and the first line of \eqref{Th1-4} can be found in the literature (see, e.g., 
Eqs.~(13.3.16)--(13.3.18) and Eqs.~(15.5.2)--(15.5.4) in \cite{NIST:DLMF};
Eqs.~(10)--(12) in page 254--255 and Eqs.~(20)--(22) in page 102 in \cite{Bateman:100233};
Eqs.~(13.4.9), (13.4.10), (13.4.13) in page 507 and Eqs.~(15.2.2)--(15.2.4) in page 557 in \cite{mabramowitz64:handbook}),
but the other identities, i.e, 
\eqref{Th1-1}, the second line of \eqref{Th1-4}, and \eqref{Th1-5}
seem to be omitted.

As we discussed in \S\ref{sec:intro}, one can formulate the identities for such exceptional cases by applying analytic continuation to the identities for the well-known cases.
However, such a derivation assumes that one already knows the identities for some cases.
Thus, this strategy may not work for more general functions, if no differentiation identities are known for those functions. 
From this point of view, we prove Theorem~\ref{thm1},
treating all the identities on an equal footing, without assuming a priori knowledge of some part of the identities.
We provide a proof of Theorem~\ref{thm1} in \S\ref{sec:proof}, which is fairly concise and is based only on the definition of the coefficients for the power series expansion of the hypergeometric functions.

In addition to the identities~\eqref{Th1-1}--\eqref{Th1-5} in Theorem~\ref{thm1}, we can generate 
other sets of identities by applying the Kummer-type transformations.
For the case of $(p,q)=(1,1)$ and $(2,1)$, we have 
\be \label{K1} \Func{1}{1}{a}{c}{z} = e^z \Func{1}{1}{c-a}{c}{-z} , \ee
and 
\begin{align}
\label{K2} \Func{2}{1}{a,b}{c}{z} &= (1-z)^{c-a-b} \Func{2}{1}{c-a,c-b}{c}{z} , \\
\label{K3} \Func{2}{1}{a,b}{c}{z} &= (1-z)^{-a} \Func{2}{1}{a,c-b}{c}{\f{z}{z-1}} ,
\end{align}
respectively, 
which enables us to obtain the following sets of identities:

\begin{landscape}
\begin{corollary} \label{coy1}
The following identities for the differentiation of the confluent hypergeometric function ${}_1F_1$ hold:
\begin{align}
&\f{d^n}{dz^n}\kk{z^r e^{-z}\Func{1}{1}{a}{c}{z} } \notag\\
&= 
\begin{cases}
      (r-n+1)_n z^{r-n} \Func{2}{2}{r+1,c-a}{r-n+1,c}{-z},  & (r\not\in \sZ_{< n}) , \\
      (-1)^{n-r} \f{n!}{(n-r)!}\f{(c-a)_{n-r}}{(c)_{n-r}} \Func{2}{2}{n+1,c-a+n-r}{n-r+1,c+n-r}{-z} , & (r=0,1,\cdots,n) , \\
      (r-n+1)_n z^{r-n} \Func{2}{2}{r+1,c-a}{r-n+1,c}{-z} + (-1)^{n-r} \f{n!}{(n-r)!}\f{(c-a)_{n-r}}{(c)_{n-r}} \Func{2}{2}{n+1,c-a+n-r}{n-r+1,c+n-r}{-z} , & (r\in \sZ_{<0}) ,
\end{cases} \label{Co1-1} \\
&\f{d^n}{dz^n} \kk{e^{-z} \Func{1}{1}{a}{c}{z} } = (-1)^n \f{(c-a)_n}{(c)_n} e^{-z} \Func{1}{1}{a}{c+n}{z} , \label{Co1-2} \\
&\f{d^n}{dz^n} \kk{z^{c-a+n-1} e^{-z} \Func{1}{1}{a}{c}{z} } = (c-a)_n z^{c-a-1} e^{-z} \Func{1}{1}{a-n}{c}{z} , \label{Co1-3} \\
&\f{d^n}{dz^n} \kk{z^{c-1} e^{-z} \Func{1}{1}{a}{c}{z} } \notag\\
&= 
\begin{cases}
(c-n)_n z^{c-n-1} e^{-z} \Func{1}{1}{a-n}{c-n}{z} , & (c\not\in \sZ_{\leq n}) , \\
(-1)^{1-c-n} \f{n!}{(n-c+1)!} \f{(c-a)_{n-c+1}}{(c)_{n-c+1}} e^{-z} \Func{1}{1}{a-c+1}{n-c+2}{z} & (c=1,2,\cdots , n+1), 
\end{cases} \label{Co1-4} \\
&\f{d^n}{dz^n} \kk{z^r e^{-z} \Func{1}{1}{c-1}{c}{z} } \notag\\
&= 
\begin{cases}
(-1)^{n-r} \f{n!}{(c)_{n-r}} e^{-z} \Func{1}{1}{c-r-1}{c-r+n}{z} & (r=1,2,\cdots , n), \\
(-1)^{n-r} \f{n!}{(c)_{n-r}} e^{-z} \Func{1}{1}{c-r-1}{c-r+n}{z} 
+ (r-n+1)_n z^{r-n} \Func{2}{2}{r+1,1}{r-n+1,c}{-z},  & (r\in \sZ_{<0}) . 
\end{cases} \label{Co1-5} 
\end{align}
\end{corollary}
\end{landscape}

\begin{landscape}
\begin{corollary} \label{coy2}
The following identities for the differentiation of the hypergeometric function ${}_2F_1$ hold:
{\fontsize{10pt}{10pt}\selectfont 
\begin{align}
&\f{d^n}{dz^n}\kk{z^r (1-z)^{a+b-c} \Func{2}{1}{a,b}{c}{z} } \notag\\
&= 
\begin{cases}
      (r-n+1)_nz^{r-n} \Func{3}{2}{r+1,c-a,c-b}{r-n+1,c}{z} , & (r\not\in \sZ_{< n}) , \\
      \f{n!}{(n-r)!}\f{(c-a)_{n-r}(c-b)_{n-r}}{(c)_{n-r}} \Func{3}{2}{n+1,c-a+n-r,c-b+n-r}{n-r+1,c+n-r}{z} , & (r=0,1,\cdots,n) , \\
      (r-n+1)_nz^{r-n} \Func{3}{2}{r+1,c-a,c-b}{r-n+1,c}{z} + \f{n!}{(n-r)!}\f{(c-a)_{n-r}(c-b)_{n-r}}{(c)_{n-r}} \Func{3}{2}{n+1,c-a+n-r,c-b+n-r}{n-r+1,c+n-r}{z} , & (r\in \sZ_{<0})
\end{cases} \label{Co2-1} \\
&\f{d^n}{dz^n}\kk{z^r (1-z)^{a+n-r-1} \Func{2}{1}{a,b}{c}{z} } \notag\\
&= 
\begin{cases}
      (r-n+1)_nz^{r-n}(1-z)^{-r-1} \Func{3}{2}{r+1,a,c-b}{r-n+1,c}{\f{z}{z-1}} , & (r\not\in \sZ_{< n}) , \\
      (-1)^{r-n}\f{n!}{(n-r)!}\f{(a)_{n-r}(c-b)_{n-r}}{(c)_{n-r}} (1-z)^{-n-1} \Func{3}{2}{n+1,a+n-r,c-b+n-r}{n-r+1,c+n-r}{\f{z}{z-1}} , & (r=0,1,\cdots,n) , \\
      (r-n+1)_nz^{r-n}(1-z)^{-r-1} \Func{3}{2}{r+1,a,c-b}{r-n+1,c}{\f{z}{z-1}} \\
      \,\,\,\,\,\, + (-1)^{r-n}\f{n!}{(n-r)!}\f{(a)_{n-r}(c-b)_{n-r}}{(c)_{n-r}} (1-z)^{-n-1} \Func{3}{2}{n+1,a+n-r,c-b+n-r}{n-r+1,c+n-r}{\f{z}{z-1}} , & (r\in \sZ_{<0})
\end{cases} \label{Co2-2} \\
&\f{d^n}{dz^n}\kk{z^{c-a+n-1}(1-z)^{a+b-c} \Func{2}{1}{a,b}{c}{z} } = (c-a)_n z^{c-a-1}(1-z)^{a+b-c-n} \Func{2}{1}{a-n,b}{c}{z} , \label{Co2-3} \\
&\f{d^n}{dz^n}\kk{ (1-z)^{a+b-c} \Func{2}{1}{a,b}{c}{z} } = \f{(c-a)_n(c-b)_n}{(c)_n} (1-z)^{a+b-c-n} \Func{2}{1}{a,b}{c+n}{z} , \label{Co2-4} \\
&\f{d^n}{dz^n}\kk{(1-z)^{a+n-1} \Func{2}{1}{a,b}{c}{z} } = (-1)^n \f{(a)_n(c-b)_n}{(c)_n} (1-z)^{a-1} \Func{2}{1}{a+n,b}{c+n}{z} , \label{Co2-5} 
\end{align}
}

{\fontsize{10pt}{10pt}\selectfont 
\begin{align}
&\f{d^n}{dz^n}\kk{z^{c-1}(1-z)^{a+b-c} \Func{2}{1}{a,b}{c}{z} } \notag\\
&= 
\begin{cases}
(c-n)_n z^{c-n-1}(1-z)^{a+b-c-n} \Func{2}{1}{a-n,b-n}{c-n}{z} , & (c\not\in \sZ_{\leq n}), \\
\f{(c-a)_{n-c+1}(c-b)_{n-c+1}}{(c)_{n-2c+2}} (1-z)^{a+b-c-n} \Func{2}{1}{a-c+1,b-c+1}{n-c+2}{z}, & (c=1,2,\cdots , n+1), \\
\end{cases} \label{Co2-6} \\
&\f{d^n}{dz^n}\kk{z^{c-1}(1-z)^{a-c+n} \Func{2}{1}{a,b}{c}{z} } \notag\\
&= 
\begin{cases}
(c-n)_n z^{c-n-1}(1-z)^{a-c} \Func{2}{1}{a,b-n}{c-n}{z} , & (c\not\in \sZ_{\leq n}), \\
(-1)^{n-c+1}\f{(a)_{n-c+1}(c-b)_{n-c+1}}{(n+1)_{1-c}(c)_{n-c+1}} (1-z)^{a-c} \Func{2}{1}{a-c+n+1,b-c+1}{n-c+2}{z}, & (c=1,2,\cdots , n+1), \\
\end{cases} \label{Co2-7} \\
&\f{d^n}{dz^n}\kk{z^r (1-z)^{a-1} \Func{2}{1}{a,c-1}{c}{z} } \notag\\
&= 
\begin{cases}
\f{n!(c-a)_{n-r}}{(c)_{n-r}} (1-z)^{a-n-1} \Func{2}{1}{a,c-r-1}{c+n-r}{z} , & (r=1,2,\cdots , n), \\
\f{n!(c-a)_{n-r}}{(c)_{n-r}} (1-z)^{a-n-1} \Func{2}{1}{a,c-r-1}{c+n-r}{z} + (r-n+1)_n z^{r-n} \Func{3}{2}{r+1,1,c-a}{r-n+1,c}{z} , & (r\in \sZ_{<0}), \\
\end{cases} \label{Co2-8} \\
&\f{d^n}{dz^n}\kk{z^r (1-z)^{a+n-r-1} \Func{2}{1}{a,c-1}{c}{z} } \notag\\
&= 
\begin{cases}
(-1)^{r-n} \f{n!(a)_{n-r}}{(c)_{n-r}} (1-z)^{a-r-1} \Func{2}{1}{a+n-r,c-r-1}{c+n-r}{z} , & (r=1,2,\cdots , n), \\
(-1)^{r-n} \f{n!(a)_{n-r}}{(c)_{n-r}} (1-z)^{a-r-1} \Func{2}{1}{a+n-r,c-r-1}{c+n-r}{z} + (r-n+1)_n z^{r-n} (1-z)^{-r-1} \Func{3}{2}{r+1,1,a}{r-n+1,c}{\f{z}{z-1}}  , & (r\in \sZ_{<0}). \\
\end{cases} \label{Co2-9} \\
&\f{d^n}{dz^n}\kk{z^r (1-z)^{n-r} \Func{2}{1}{1,a}{c}{z} } \notag\\
&= 
\begin{cases}
(-1)^{r-n} \f{n!(c-a)_{n-r}}{(c)_{n-r}} \Func{2}{1}{n+1,a}{c+n-r}{z} , & (r=1,2,\cdots , n), \\
(-1)^{r-n} \f{n!(c-a)_{n-r}}{(c)_{n-r}} \Func{2}{1}{n+1,a}{c+n-r}{z} + (r-n+1)_n z^{r-n} (1-z)^{-r-1} \Func{3}{2}{r+1,1,c-a}{r-n+1,c}{\f{z}{z-1}} , & (r\in \sZ_{<0}). \\
\end{cases} \label{Co2-10}
\end{align}
}
\end{corollary}

\end{landscape}

Again, some of these identities are not written down explicitly in the literature.
For Corollary~\ref{coy1}, \eqref{Co1-2} and \eqref{Co1-3} 
correspond to Eqs.~(13.3.19) and (13.3.20) in \cite{NIST:DLMF}, 
Eqs.~(13) and (14) in page 255 of \cite{Bateman:100233}, or 
Eqs.~(13.4.12) and (13.4.11) in page 507 in \cite{mabramowitz64:handbook}.
An identity corresponding to the first line of \eqref{Co1-4} can be found in Eq.~(13.3.21) in \cite{NIST:DLMF} or Eq.~(13.4.14) in page 507 in \cite{mabramowitz64:handbook}. 
For Corollary~\ref{coy2}, \eqref{Co2-3}--\eqref{Co2-5}, the first line of \eqref{Co2-6}, and the first line of \eqref{Co2-7} correspond to Eqs.~(15.5.5)--(15.5.9) in \cite{NIST:DLMF}, Eqs.~(23)--(27) in page 102--103 in \cite{Bateman:100233}, or Eqs.~(15.2.5)--(15.2.9) in page 557 in \cite{mabramowitz64:handbook}. 
However, the other identities in Corollaries~\ref{coy1} and \ref{coy2}, i.e., \eqref{Co1-1}, the second line of \eqref{Co1-4}, \eqref{Co1-5}--\eqref{Co2-2}, the second line of \eqref{Co2-6}, the second line of \eqref{Co2-7}, and \eqref{Co2-8}--\eqref{Co2-10} are not explicitly mentioned in the literature.

The proof of Corollaries~\ref{coy1} and \ref{coy2} is straightforward.
Note that, after applying the Kummer-type transformations~\eqref{K1}--\eqref{K3} to the identities~\eqref{Th1-1}--\eqref{Th1-5} in Theorem~\ref{thm1}, one needs to rename the variable and parameters appropriately.
Specifically, after applying \eqref{K1}, one should perform the following replacements: $a\to c-a$ and $z\to -z$. 
Similarly, after applying \eqref{K2}, the necessary replacements are $a\to c-a$ and $b\to c-b$.
Finally, after applying \eqref{K3}, one should rename $b\to c-b$ and $z\to \f{z}{z-1}$, and hence
\be \f{d^n}{dz^n} \to (-1)^n (1-z)^{n+1} \f{d^n}{dz^n} (1-z)^{n-1} . \ee

More generally, one can also apply
generalized Kummer-type transformations for ${}_{p+1}F_{p}$ \cite{Miller2011} and ${}_{p}F_{p}$ \cite{MILLER2009964}
to generate similar identities.

\section{Proof of Theorem \ref{thm1}}
\label{sec:proof}

We consider a differentiation 
\be \label{diff} \f{d^n}{dz^n}\kk{z^r \Func{p}{q}{\ba}{\bb}{z} } = \f{d^n}{dz^n} \sum_{k=0}^\infty c_k z^{k+r}. \ee
According to which terms in the series survive after performing the term-by-term differentiation in the right-hand side, we consider the following three cases separately:
\begin{enumerate}
\item $r\not\in \sZ_{<n}$,
\item $r=0,1,\cdots, n$,
\item $r\in \sZ_{<0}$.
\end{enumerate}
Note that the first two cases overlap for $r=n$, corresponding to the similar overlapping in \eqref{Th1-1} and \eqref{Th1-4}, for which the formulae coincide.
Below we prove the identity~\eqref{Th1-1} for each case, from which we derive other identities~\eqref{Th1-2}--\eqref{Th1-5} by using the reduction~\eqref{cancel}.

\subsection{$r\not\in \sZ_{<n}$}
\label{ssec:1}

In this case, all the terms in the series in the right-hand side of \eqref{diff} survive after performing the term-by-term differentiation. 
We thus obtain 
\begin{align}
\f{d^n}{dz^n} z^{k+r} &= (k+r)(k+r-1)\cdots (k+r-n+1)z^{k+r-n} \notag\\
&= (k+r-n+1)_n z^{k+r-n} .
\end{align}
Using the definitions~\eqref{pFqdef} and \eqref{ckdef}, and an identity
\be \label{id1} \f{(a+k)_n}{(a)_n} = \f{(a+n)_k}{(a)_k} , \ee
we obtain
\begin{align} 
\f{d^n}{dz^n}\kk{z^r \Func{p}{q}{\ba}{\bb}{z} } &= \sum_{k=0}^\infty \f{1}{k!}\f{(\ba)_k}{(\bb)_k} (k+r-n+1)_n z^{k+r-n} \notag\\ 
&= (r-n+1)_n z^{r-n} \sum_{k=0}^\infty \f{1}{k!} \f{(r+1,\ba)_k}{(r-n+1,\bb)_k} z^k \notag\\
&= (r-n+1)_nz^{r-n} \Func{p+1}{q+1}{r+1,\ba}{r-n+1,\bb}{z},
\label{case1} 
\end{align}
which establishes the first line of \eqref{Th1-1} for the case $r\not\in \sZ_{<n}$.

By virtue of the relation~\eqref{cancel}, $\Func{p+1}{q+1}{r+1,\ba}{r-n+1,\bb}{z}$ in the right-hand side of \eqref{case1} can be reduced to $\Func{p}{q}{\ba'}{\bb'}{z}$ with some $\ba'$ and $\bb'$ if an upper parameter coincides with a lower parameter.
While such ``cancellation'' between argument(s) of $\ba$ and $\bb$ may be possible, it should be studied on a case-by-case basis.
Here, we focus on nontrivial cases peculiar to \eqref{case1}, which show up by a cancellation involving $r+1$ or $r-n+1$.
Without loss of generality, we assume that the first argument of $\ba$ or $\bb$ is a counterpart of the cancellation.

Clearly, there are two ways of cancellation.
If $r=a_1+n-1$, $a_1$ and $r-n+1$ cancel each other, and 
then \eqref{case1} reduces to the identity~\eqref{Th1-3}.
Note that, however, the condition $r\not\in \sZ_{<n}$ now reads $a_1\not\in \sZ_{\leq 0}$, so \eqref{Th1-3} for $a_1\in \sZ_{\leq 0}$ needs to be proved separately.
For the case $a_1=0$, \eqref{Th1-3} holds trivially. 
We shall provide the proof of \eqref{Th1-3} for the remaining case~$a_1\in \sZ_{< 0}$ in \S\ref{ssec:3}.

On the other hand, if $r=b_1-1$, $b_1$ and $r+1$ cancel each other, and 
then \eqref{case1} reduces to the first line of \eqref{Th1-4}.
The condition $r\not\in \sZ_{<n}$ prohibits $b_1-n$ to be a nonpositive integer, 
which is implicitly assumed for the right-hand side of the first line of \eqref{Th1-4} to avoid the singularity of the coefficients.

\subsection{$r=0,1,\cdots, n$}
\label{ssec:2}

In this case, the first $n-r$ terms vanish after performing the term-by-term differentiation in \eqref{diff}. 
As a result, we obtain
\begin{align} 
\f{d^n}{dz^n}\kk{z^r \Func{p}{q}{\ba}{\bb}{z} } 
&= n!c_{n-r} \sum_{k=0}^\infty \f{(n+k)!}{n!k!} \f{c_{n-r+k}}{c_{n-r}} z^k \notag\\
&= n!c_{n-r} \sum_{k=0}^\infty \f{(n+1)_k}{k!} \f{(n-r)!}{(n-r+k)!}\f{(\ba)_{n-r+k}}{(\ba)_{n-r}}\f{(\bb)_{n-r}}{(\bb)_{n-r+k}} z^k \notag\\
&= \f{n!}{(n-r)!} \f{(\ba)_{n-r}}{(\bb)_{n-r}} \sum_{k=0}^\infty \f{1}{k!} \f{(n+1,\ba+n-r)_k}{(n-r+1,\bb+n-r)_k} z^k \notag\\
&= \f{n!}{(n-r)!}\f{(\ba)_{n-r}}{(\bb)_{n-r}} \Func{p+1}{q+1}{n+1,\ba+n-r}{n-r+1,\bb+n-r}{z}, 
\label{case2}
\end{align}
where we used an identity
\be \label{id2} \f{(a)_{k+n-r}}{(a)_{n-r}} = (a+n-r)_{k}. \ee
The right-hand side of \eqref{case2} is nothing but the second line of \eqref{Th1-1}.

Again, let us consider nontrivial reduction of ${}_{p+1}F_{q+1}$ in the right-hand side of \eqref{case2} to ${}_{p}F_{q}$ by using \eqref{cancel}.
There are three ways of nontrivial cancellation.
The simplest cancellation occurs if $r=0$, for which $n+1$ and $n-r+1$ cancel.
In this case, after cancellation, \eqref{case2} yields \eqref{Th1-2}.

For $r\ne 0$, in parallel to the case considered in \S\ref{ssec:1}, the factor $n+1$ or $n-r+1$ cancel with $b_1+n-r$ or $a_1+n-r$ respectively, where we focused on the first arguments of $\ba$ and $\bb$ without loss of generality.
If $r=b_1-1$, where $b_1=1,2,\cdots,n+1$ is required, 
\eqref{case2} yields the second line of \eqref{Th1-4}.
On the other hand, if $a_1=1$, \eqref{case2} yields \eqref{Th1-5}.

\subsection{$r\in \sZ_{<0}$}
\label{ssec:3}

Let us denote $r=-m\,(m\in \sN)$ in this case.
We can decompose the series in \eqref{diff} into three parts:
\begin{align} 
\label{case3} \f{d^n}{dz^n}\kk{z^{-m} \Func{p}{q}{\ba}{\bb}{z} } 
= \f{d^n}{dz^n} ( &c_0 z^{-m} + c_1 z + \cdots + c_{m-1} z^{-1} \notag\\
&+ c_m + c_{m+1}z + \cdots + c_{m+n-1} z^{n-1} \notag\\
&+ c_{m+n}z^n + c_{m+n+1}z^{n+1} + \cdots ) .
\end{align}
Terms in the second line of the right-hand side of \eqref{case3} vanish after performing the term-by-term differentiation, whereas the first and third lines survive in general.
As a result, we obtain
\begin{align} 
\label{case3-2pre} 
\f{d^n}{dz^n}\kk{z^{-m} \Func{p}{q}{\ba}{\bb}{z} } 
&= (-m-n+1)_n z^{-m-n} \sum_{k=0}^{m-1} \f{(-m-n+1+k)_n}{(-m-n+1)_n} c_k z^k \notag\\
&\,\,\,\,\,\, + n! c_{m+n} \sum_{k=0}^\infty \f{(k+1)_n}{(1)_n} \f{c_{m+n+k}}{c_{m+n}} z^k ,
\end{align}
where we normalized the coefficients of the series to $1$ for $k=0$.
Using the definition~\eqref{ckdef} and the identity~\eqref{id1}, the coefficients simplifies as 
\begin{align} 
\label{coe1} \f{(-m-n+1+k)_n}{(-m-n+1)_n} c_k &= \f{1}{k!} \f{(-m+1,\ba)_k}{(-m-n+1,\bb)_k}, \\
\label{coe2} \f{(k+1)_n}{(1)_n} \f{c_{m+n+k}}{c_{m+n}} &= \f{1}{k!} \f{(n+1,\ba+m+n)_k}{(m+n+1,\bb+m+n)_k} .
\end{align}
Plugging these expressions, we see that \eqref{case3-2pre} yields the third line of \eqref{Th1-1}, namely,
\begin{align} 
\label{case3-2} 
\f{d^n}{dz^n}\kk{z^{-m} \Func{p}{q}{\ba}{\bb}{z} } 
&= (-m-n+1)_nz^{-m-n} \Func{p+1}{q+1}{-m+1,\ba}{-m-n+1,\bb}{z} \notag\\
&\,\,\,\,\,\, + \f{n!}{(n+m)!}\f{(\ba)_{n+m}}{(\bb)_{n+m}} \Func{p+1}{q+1}{n+1,\ba+n+m}{n+m+1,\bb+n+m}{z} .
\end{align}
Note that, while the first series in the right-hand side of \eqref{case3-2pre} is a finite series up to the order of $z^{m-1}$, it precisely yields the first generalized hypergeometric function in the right-hand side of \eqref{case3-2}. 
This is because $\Func{p+1}{q+1}{-m+1,\ba}{-m-n+1,\bb}{z}$ contains the nonpositive integer $-m+1$ in the upper row, and hence, by definition~\eqref{inter}, it is $(m-1)$-th order polynomial.

Let us consider a possibility of nontrivial reductions of ${}_{p+1}F_{q+1}$ in the right-hand side of \eqref{case3-2} to ${}_{p}F_{q}$ by using \eqref{cancel}.
Again, we focus on the first arguments of $\ba$, $\bb$, $\ba+n+m$ or $\bb+n+m$ without loss of generality.
First, for the first term in the right-hand side of \eqref{case3-2}, 
in general it is impossible to cancel $-m+1$ with $b_1$ 
(unless $\ba$ contains nonpositive integer(s)). 
The remaining possibility of nontrivial cancellation is the one between $-m-n+1$ and $a_1$, which occurs if $a_1=-m-n+1$.
In this case, the original $\Func{p}{q}{\ba}{\bb}{z}$ in the left-hand side is $(m+n-1)$-th order polynomial.
Hence, the third line of the right-hand side of \eqref{case3} vanishes after performing the term-by-term differentiation, which implies that the second hypergeometric function in the right-hand side of \eqref{case3-2} does not exist.
Consequently, we obtain 
\be \f{d^n}{dz^n}\kk{z^{-m} \Func{p}{q}{\ba}{\bb}{z} } = (-m-n+1)_n z^{-m-n} \Func{p}{q}{-m+1,\hat \ba}{\bb}{z}, \ee
which, after substituting $m=-a_1-n+1$, yields \eqref{Th1-3} for the case $a_1\in \sZ_{<0}$.

Finally, we consider a cancellation in the second term in the right-hand side of \eqref{case3-2}.
Again, the cancellation between $n+1$ and $b_1+m+n$ is impossible in general to avoid the singular behavior of, in this case, the left-hand side. 
The cancellation between $m+n+1$ and $a_1+m+n$ occurs if $a_1=1$.
In this case, \eqref{case3-2} yields the second line of \eqref{Th1-5}. 
\qed

\section*{Acknowledgments}

This work was supported by Japan Society for the Promotion of Science (JSPS) Grants-in-Aid for Scientific Research (KAKENHI) Grant No.~JP18K13565 and No.~JP22K03639.

\bibliographystyle{elsarticle-num}
\bibliography{refs}

\begin{thebibliography}{10}
\expandafter\ifx\csname url\endcsname\relax
  \def\url#1{\texttt{#1}}\fi
\expandafter\ifx\csname urlprefix\endcsname\relax\def\urlprefix{URL }\fi
\expandafter\ifx\csname href\endcsname\relax
  \def\href#1#2{#2} \def\path#1{#1}\fi

\bibitem{Bailey}
W.~N. Bailey, {Generalized Hypergeometric Series}, Cambridge Tracts in
  Mathematics, Cambridge Univ. Press, Cambridge, 1935.

\bibitem{Slater}
L.~J. Slater, {Generalized Hypergeometric Functions}, 2nd Edition, Cambridge
  Univ. Press, Cambridge, 1966.

\bibitem{Bateman:100233}
H.~Bateman, A.~Erdélyi,
  \href{https://resolver.caltech.edu/CaltechAUTHORS:20140123-104529738}{{Higher
  transcendental functions, Volume I}}, California Institute of technology.
  Bateman Manuscript project, McGraw-Hill, New York, NY, 1955.
\newline\urlprefix\url{https://resolver.caltech.edu/CaltechAUTHORS:20140123-104529738}

\bibitem{mabramowitz64:handbook}
M.~Abramowitz, I.~A. Stegun (Eds.), Handbook of Mathematical Functions with
  Formulas, Graphs and Mathematical Tables, Dover Publications, Inc., New York,
  1965.

\bibitem{Luke}
Y.~Luke, {The Special Functions and Their Approximations, Vol. 1}, Academic
  Press, 1969.

\bibitem{NIST:DLMF}
\href{http://dlmf.nist.gov/}{{\it NIST Digital Library of Mathematical
  Functions}}, Release 1.1.3 of 2021-09-15, F.~W.~J. Olver, A.~B. {Olde
  Daalhuis}, D.~W. Lozier, B.~I. Schneider, R.~F. Boisvert, C.~W. Clark, B.~R.
  Miller, B.~V. Saunders, H.~S. Cohl, and M.~A. McClain, eds.
\newline\urlprefix\url{http://dlmf.nist.gov/}

\bibitem{Seaborn}
J.~B. Seaborn, {Hypergeometric Functions and Their Applications}, Texts in
  Applied Mathematics, Springer-Verlag, 1991.

\bibitem{Motohashi:prep}
H.~Motohashi, {Differentiation identities for Heun functions}, in preparation.

\bibitem{Magnus}
W.~Magnus, F.~Oberhettinger, R.~P. Soni, {Formulas and Theorems for the Special
  Functions of Mathematical Physics}, Die Grundlehren der mathematischen
  Wissenschaften, Springer-Verlag, 1966.

\bibitem{Miller2011}
A.~R. Miller, R.~B. Paris, Euler-type transformations for the generalized
  hypergeometric function ${}_{r+2}F_{r+1}(x)$, Z. Angew. Math. Phys. 62 (2011)
  31--45.
\newblock \href {http://dx.doi.org/https://doi.org/10.1007/s00033-010-0085-0}
  {\path{doi:https://doi.org/10.1007/s00033-010-0085-0}}.

\bibitem{MILLER2009964}
A.~R. Miller,
  \href{https://www.sciencedirect.com/science/article/pii/S0377042709003252}{Certain
  summation and transformation formulas for generalized hypergeometric series},
  Journal of Computational and Applied Mathematics 231~(2) (2009) 964--972.
\newblock \href {http://dx.doi.org/https://doi.org/10.1016/j.cam.2009.05.013}
  {\path{doi:https://doi.org/10.1016/j.cam.2009.05.013}}.
\newline\urlprefix\url{https://www.sciencedirect.com/science/article/pii/S0377042709003252}

\end{thebibliography}

\end{document}